\documentclass[12pt, reqno]{amsart}
\usepackage[dvips]{graphicx}
\usepackage{amssymb,latexsym, amsmath, amscd, array}

\numberwithin{equation}{section} 
\numberwithin{figure}{section}
\numberwithin{table}{section}

\DeclareMathOperator{\vol}{{\rm vol}}
\DeclareMathOperator{\iqvol}{{\rm IQvol}}

\DeclareMathOperator{\area}{{\rm area}}

\DeclareMathOperator{\stsys}{{\rm stsys}}
\DeclareMathOperator{\sys}{{\rm sys}} 

\DeclareMathOperator{\sysh}{{\rm sysh}}

\DeclareMathOperator{\cuplength}{{\rm cup-length}}

\DeclareMathOperator{\syscat}{{\rm cat_{sys}}}

\DeclareMathOperator{\IQ}{{\rm IQ}}
\DeclareMathOperator{\cat}{{\mbox{\rm cat$_{\rm LS}$}}}
\DeclareMathOperator{\dR}{{\rm dR}}
\DeclareMathOperator{\Ker}{{\rm Ker}}
\DeclareMathOperator{\Hom}{{\rm Hom}}
\DeclareMathOperator{\IM}{{\rm Im}}

\DeclareMathOperator{\In}{{\rm Indet}}
\DeclareMathOperator{\AJ} {{\mathcal A}}

\DeclareMathOperator{\pisys}{{\rm sys}\pi}

\def\gmetric{{\mathcal G}}

\def\rp{\mathbb R\mathbb P}
\def\cp{\mathbb C\mathbb P}

\def\la{\langle} 
\def\ra{\rangle} 
\def\ov{\overline} 

\def\FF{{F_M}}

\newcommand{\C}{\mathbb C} 
\newcommand{\N}{\mathbb N}

\newcommand{\R}{\mathbb R} 

\newcommand{\Z}{\mathbb Z}

\newtheorem{theorem}{Theorem}[section]

\newtheorem{lemma}[theorem]{Lemma}
\newtheorem{corollary}[theorem]{Corollary}

\newtheorem{prop}[theorem]{Proposition}

\newtheorem{cory}[theorem]{Corollary}

\theoremstyle{definition} 
\newtheorem{definition}[theorem]{Definition}
\newtheorem{example}[theorem]{Example}
\newtheorem{remark}[theorem]{Remark}
\newtheorem{rem}[theorem]{Remark}
\newtheorem{df}[theorem]{Definition}
\newtheorem{question}[theorem]{Question} 

\newcommand\theoref{Theorem~\ref} 
\newcommand\lemref{Lemma~\ref}
\newcommand\propref{Proposition~\ref}
\newcommand\coryref{Corollary~\ref} 

\def\ie {{\it i.e.\ }} \def\etc {{\it etc}} \def\eg {{\it e.g.\ }}
\def\cf {\hbox{\it cf.\ }}

\def\ga{\alpha}

\def\eps{\varepsilon} 
\def\gf{\varphi}

\def\m{\medskip}

\date{\today}
\begin{document}

\author[M.~Katz]{Mikhail G. Katz$^{*}$} \address{Department of
Mathematics, Bar Ilan University, Ramat Gan 52900 Israel}
\email{katzmik@math.biu.ac.il} \thanks{$^{*}$Supported by the Israel
Science Foundation (grants no.\ 620/00-10.0 and 84/03)}

\author[Y.~Rudyak]{Yuli B. Rudyak$^{\dagger}$} \address{Department of
Mathematics, University of Florida, PO Box 118105, Gainesville, FL
32611-8105 USA} \email{rudyak@math.ufl.edu}
\thanks{$^{\dagger}$Supported by MCyT, project BFM 2002-00788, Spain, and by 
NSF, grant 0406311}

\title[Lusternik-Schnirelmann meet systole]{Lusternik-Schnirelmann
category and systolic category of low dimensional
manifolds$^1$}

\subjclass
{Primary 53C23;  
Secondary  55M30, 
57N65  
}

\keywords{detecting element, essential manifolds, isoperimetric
quotient, Lusternik-Schnirelmann category, Massey product, systole}

\begin{abstract}
We show that the geometry of a Riemannian manifold~$(M,\gmetric)$ is
sensitive to the apparently purely homotopy-theoretic invariant of $M$
known as the Lusternik-Schnirelmann category, denoted $\cat (M)$.
Here we introduce a Riemannian analogue of~$\cat (M)$, called the {\em
systolic category\/} of~$M$.  It is denoted $\syscat (M)$, and defined
in terms of the existence of systolic inequalities satisfied by every
metric $\gmetric$, as initiated by C.~Loewner and later developed by
M.~Gromov.  We compare the two categories.  In all our examples, the
inequality $\syscat M \le \cat M$ is satisfied, which typically turns
out to be an equality, \eg in dimension 3.  We show that a number of
existing systolic inequalities can be reinterpreted as special cases
of such equality, and that both categories are sensitive to Massey
products.  The comparison with the value of~$\cat (M)$ leads us to
prove or conjecture new systolic inequalities on $M$.
\end{abstract}

\footnotetext [1]{{\em Communications on Pure and Applied
Mathematics}, to appear.  Available at \texttt{arXiv:math.DG/0410456}}

\maketitle

\tableofcontents

\section*{Introduction}

In his paper \cite{G}, T. Ganea demonstrated the usefulness of the
study of numerical (\ie $\Z$- or $\N$-valued) topological invariants, 
such as the Lusternik--Schnirelmann category
$\cat$. For closed smooth manifolds $M$, one can consider additional
numerical invariants, such as the minimal number of balls that cover
the manifold, or the minimal number of critical points of a smooth
function on $M$, \etc.  Here we introduce (a few versions of) a
Riemannian analogue of $\cat(M)$, called the systolic category of $M$,
as in Definition~\ref{syscat}.  It is defined in terms of the
existence of ``curvature-free'' systolic inequalities satisfied by
every metric $\gmetric$, as initiated by C.~Loewner~\cite{Pu} and
developed later by M.~Gromov \cite{Gr1}.

\m Given that a dozen or so numerical invariants have by now been studied
\cite{CLOT}, one could legitimately ask, why define another one?  We
feel that systolic category is the only numerical invariant possessing a
differential geometric flavor, thus adding some Riemannian spice to an
increasingly homotopy-theoretic field.  Furthermore, systolic category
provides an intuitive point of entry for differential geometers into
the field of numerical invariants.  

\m It is natural to compare $\syscat$ with the other numerical
invariants mentioned above.  Here we start this program and show that
$\syscat$ is sensitive to $\cat$, which is a purely homotopic
invariant.  We verify the equality $\syscat=\cat$ in dimension 3 based
on the results of G\'omez-Larra\~naga and Gonz\'ales-Acu\~na
\cite{GG}, see also \cite{OR}.  We show that a number of existing
systolic inequalities can be reinterpreted as special cases of such
equality, \cf \eqref{33}, \eqref{52}, and \eqref{85}.

We state at the outset that, while a precise Definition~\ref{syscat}
appears below, it is not entirely clear what the ``right'' definition
of systolic category should be, \cf \eqref{52}, \eqref{33},
\eqref{94}.  We present what we feel is compelling evidence suggesting
the existence of a coherent notion of this sort.

In higher dimensions, no immediate (direct) connections of the two
categories are as yet available, but we observe that they have
parallel behavior and common points of sensitivity.  More
specifically, we make the following observations:

\begin{itemize}
\item
the two categories coincide in dimensions 2 and 3;
\item
the two categories attain the maximal value (\ie dimension)
simultaneously;
\item
both categories admit a lower bound in terms of the real cup-length;
\item
both categories are sensitive to Massey products, \cf
section~\ref{massey}.
\end{itemize}

The comparison with the value of~$\cat (M)$ leads us to prove or
conjecture new systolic inequalities on~$M$, \eg \eqref{71y},
\eqref{71}, \eqref{101}, \eqref{86}.  Conversely, the comparison with
the value of $\syscat$ leads to open questions about $\cat$, \cf
\eqref{73}.

\m All manifolds are assumed to be closed, connected and smooth unless
explicitly mentioned otherwise.  To the extent that our paper aims to
address both a topological and a geometric audience, we attempt to
give some indication of proof of pertinent results that may be more
familiar to one audience than the other.

\section{Systoles}

Let $M$ be a (smooth) manifold equipped with a Riemannian
metric~$\gmetric$.  We will now define the systolic invariants
$\sys_k(M,\gmetric)$.  In the sequel, we abbreviate
$\sys_k(M,\gmetric)$ as $\sys_k(M)$, $\sys_k(\gmetric)$ or even
$\sys_k$, depending on the context.

\begin{definition}
The homotopy 1-systole, denoted $\pisys_1(M,\gmetric)$, is the least
length of a non-contractible loop in $M$.  The homology 1-systole,
denoted $\sysh_1(M,\gmetric)$, is defined in a similar way, in terms
of loops which are not zero-homologous.
\end{definition}

Clearly, $\pisys_1\le \sysh_1$.  Now let $k\in \N$.  Higher homology
$k$-systoles $\sysh_k$, with coefficients over a ring $A=\Z$ or
$\Z_2$, are defined similarly to $\sysh_1$, as the infimum of
$k$-areas of~$k$-cycles, with coefficients in $A$, which are not
zero-homologous.  Note that we adopt the usual convention, convenient
for our purposes, that the infimum over an empty set is infinity.  

When $k=n$ is the dimension, $\sysh_n(M,\gmetric)$ is equal to the
total volume $\vol_n(M,\gmetric)$ of a compact Riemannian $n$-manifold
$(M,\gmetric)$.  More detailed definitions appear in the survey
\cite{CK}.  We do not consider higher ``homotopy'' systoles.  All
systolic notions can be defined similarly on polyhedra, as
well~\cite{Gr2, Ba02}.

\begin{definition}[\cf \cite{Fe, BanK}]
Given a class $\ga\in H_k(M;\Z)$ of infinite order, we define the
stable norm $\| \ga_\R \|$ by setting
$$
\|\ga_\R\|=\lim_{m\to \infty} m^{-1} \inf_{\ga(m)} \vol_k(\ga(m)),
$$
where $\alpha_\R$ denotes the image of $\alpha$ in real homology,
while $\ga(m)$ runs over all Lipschitz cycles with integral
coefficients representing $m\ga$. The stable homology $k$-systole,
denoted $\stsys_k(\gmetric)$, is defined by minimizing the stable norm
$\| \alpha_\R \|$ over all integral $k$-homology classes $\alpha$ of
infinite order.
\end{definition}

\begin{remark}
All the systolic invariants defined in this section are positive for
smooth compact manifolds.  For the homology and homotopy 1-systoles,
this follows by comparison with the injectivity radius.  For stable
$k$-systoles the positivity follows by a calibration argument.  For
the $k$-systoles ($k\geq 2$), since torsion classes are involved, one
cannot use differential forms.  An argument proving positivity was
outlined in \cite[Section~2]{BCIK2}.
\end{remark}

The stable homology 1-systole $\stsys_1$ is easy to compute for
orientable surfaces, because it coincides with the ordinary homology
1-systole.  In fact this is true for the codimension one systole of
any orientable manifold \cite{Fe2, Wh}.  However, the situation is
entirely different when the codimension is bigger than 1.  For
example, every 3-manifold~$M$ with $b_1(M)\geq 1$ is
(1,2)-systolically free \cite{BabK}, meaning that there exists a
sequence of metrics $\gmetric_j$ of volume going to zero as $j
\to\infty$, and yet the product of systoles is big, \ie
\begin{equation}
\label{11}
\sysh_1(\gmetric_j) \sysh_2(\gmetric_j) > 1.
\end{equation}
Meanwhile, the freedom phenomenon disappears when we replace $\sysh_1$
by the stable systole $\stsys_1$, \cf inequality~\eqref{gromov}.

\section{Systolic categories}

Let $(M^n,\gmetric)$ be a Riemannian manifold.  Recall that, in our
convention, the systolic invariants are infinite when defined over an
empty set of loops or cycles.

\begin{definition}
\label{sys}
Given $k\in \N, k>1$ we set
$$
\sys_{k}(\gmetric)=\min\{\sysh_k(\gmetric,\Z),
\sysh_k(\gmetric,\Z_2),\stsys_k(\gmetric)\} .
$$
Furthermore, we define
$$
\sys_{1}(\gmetric)=\min\{\pisys_1(\gmetric),\sysh_1(\gmetric,\Z),
\sysh_1(\gmetric,\Z_2),\stsys_1(\gmetric)\}.
$$
\end{definition}

Let~$d\geq 2$ be an integer.  Consider a partition
\begin{equation}
\label{2.1}
n= k_1 + \ldots + k_d,
\end{equation}
where $k_i\geq 1$ for all $i=1,\ldots, d$.  We will consider
scale-invariant inequalities ``of length $d$'' of the following type:
\begin{equation}
\label{dd}
\sys_{k_1}(\gmetric) \sys_{k_2}(\gmetric) \ldots \sys_{k_d}(\gmetric) \leq
C(M) \vol_n(\gmetric),
\end{equation}
satisfied by all metrics $\gmetric$ on $M$, where the constant $C(M)$
depends only on the topological type of $M$ but not on the metric
$\gmetric$.  Here $\sys_k$ denotes the minimum of all non-vanishing
systolic invariants in dimension $k$ as defined above. 

\begin{definition}
\label{syscat}
Systolic category of $M$, denoted $\syscat(M)$, is the largest
integer~$d$ such that there exists a partition~\eqref{2.1} with
\[
\prod\limits^{d}_{i=1} \sys_{k_i}(M,\gmetric) \leq C(M) \vol_n(M,\gmetric) 
\]
for all metrics $\gmetric$ on $M$.
\end{definition}

In particular, $\syscat M \le \dim M$.
\begin{example}
Every orientable $n$-manifold $N$ with positive Betti number
$b_k(N)\geq 1$ satisfies Gromov's stable systolic inequality
\begin{equation}
\label{gromov}
\stsys_k(N) \stsys_{n-k}(N) < C(b_k(N)) \vol_n (N),
\end{equation}
In our notation, inequality~\eqref{gromov} implies the following bound
for the stable systolic category:
\begin{equation}
\label{52}
\syscat(N) \geq 2 \hbox{\ \ if\ } b_k(N) \geq 1 \hbox{\ for\ some\ }
1\leq k \leq n-1,
\end{equation}
see \eqref{85} for a generalisation in terms of real cuplength.  The
constant in~\eqref{gromov} depends on the Betti number, and was
studied in \cite{He, BanK}.  The constant can be optimal in certain
cases, such as Gromov's inequality for the stable 2-systole of complex
projective spaces \eqref{61}, as well as for $k=1$ \cite{BanK, BaKa}.
\end{example}

In general, systolic inequalities involving any higher $k$-systole are
notoriously hard to come by, in view of the widespread phenomenon of
systolic freedom, see \eqref{11}, \cite[Systolic reminiscences,
p.~271]{Gr3}, as well as \cite[Appendix D]{Gr3}, \cite{Ba02, Ka}.

\begin{example}
It is still unknown whether $\rp^3$ admits a (1,2)-systolic
inequality, or whether on the contrary it is (1,2)-systolically free.
Here, of course, the systoles are over $\Z_2$.  The analogous question
for the 3-manifold~$S^1\times S^2$ was resolved, in favor of freedom,
by M.~Freedman~\cite{Fr}.  The case of $\rp^3$ is not interesting from
the vantage point of category, since it is essential and therefore
\[
\cat(\rp^3) =\syscat(\rp^3)=3
\]
and not 2, see \theoref{equiv}.  Thus, the question of the existence
of a~(1,2)-systolic inequality does not affect the value of the
categories.  On the other hand, the case of $\rp^2 \times S^2$ is a
simple case that's completely open, \cf Example~\ref{112}.  Smale's
spin rational homology 5-sphere $M_k$ (see \cite{Sm}) of
Example~\ref{smale} satisfies $\cat(M_k)= 2$ by Theorem \ref{condim},
but its systolic category is inaccessible with the techniques
available.
\end{example}

We express the hope that by exhibiting a connection with
Lusternik-Schnirelmann category as well as some motivated conjectures,
we will stimulate further research in the direction of new systolic
inequalities.  We therefore conclude with the following questions.

\begin{question}
To what extent can Gromov's filling techniques~\cite{Gr1} be
generalized in the direction of proving systolic inequalities
involving higher $k$-systoles?
\end{question}

\begin{question}
It is clear that $\syscat (M\times N)\le \syscat(M) + \syscat(N)$. How
far can the inequality be from equality?  Are there examples with
$\syscat (M\times S^k)=\syscat M$?
\end{question}

\section{Categories agree in dimension 2}

Every compact surface $M$ (orientable or not) with infinite
fundamental group satisfies Gromov's inequality \cite{Gr1}
\[
\pisys_1(M)^2 \leq \frac{4}{3} \area(M).
\]
For the projective plane $\rp^2$, we have Pu's optimal inequality
\cite{Pu}
\begin{equation}
\label{21}
(\pisys_1(\rp^2))^2 \leq \frac{\pi}{2} \area(\rp^2).
\end{equation}
Tighter estimates are available as Euler characteristic becomes unbounded
\cite{Gr1, KS2}, but they are irrelevant here.  The conclusion is that
\begin{equation}
\label{amazingconjecture}
\cat M = \syscat M
\end{equation}
for all closed surfaces $M$.  Namely, we have $\cat M= 1$ for $M=S^2$
and $\cat M=2$ for all~$M\not= S^2$.  However, one cannot expect
\eqref{amazingconjecture} to hold in all dimensions, \cf
Example~\ref{singhof}.

\section{Essential manifolds and detecting elements}

A closed $n$-manifold $M$ is called {\em essential\/} if there exist a
group~$\pi$, a (possibly twisted) coefficient system $A$ on the
Eilenberg--Mac\!~Lane space $K(\pi,1)$, and a map $f: M \to K(\pi,1)$,
such the homomorphism
\[
f^*: H^n(K(\pi,1);A) \to H^n(M;f^*(A))
\]
is non-zero. Here without loss of generality we can assume that
$\pi=\pi_1(M)$. When we want to fix the twisted system $A$, we say
that $M$ is $A$-essential. In the language of \cite{OR,R}, $M$ is
$A$-essential if and only if the fundamental cohomology class in
$H^n(M;A)$ is a detecting element of category weight $n$. Clearly, $M$
is essential if it admits a map to a space $K(\pi,1)$ such
that the induced homomorphism in $n$-dimensional homology sends the
$A$-fundamental class to a nonzero class.

\begin{theorem}\label{equiv}
Let $M$ be a manifold of dimension $n$. The following three conditions
are equivalent:
\begin{enumerate}
\item[(1)] the manifold $M$ is essential;
\item[(2)] we have $\cat M=n=\dim M$;
\item[(3)] if a map $f: M \to K(\pi_1(M),1)$ induces an isomorphism of
fundamental groups, then $f(M)$ is not contained in the
$(n-1)$-skeleton of $K(\pi_1(M),1)$.
\end{enumerate}
\end{theorem}

\begin{proof}
The equivalence of (1) and (2) is proved in \cite{Ber},
see also \cite[Theorem~2.51]{CLOT}. In fact, I. Berstein \cite{Ber} proved
that $\cat M= \dim M$ if and only if $M$ is essential in a particular
sheaf, namely the tensor product $I \otimes \ldots \otimes I$ of the
augmentation ideal $I$ of the group ring of the fundamental group.  In
greater detail, one has an element $a\in H^1(M;I)$, the characteristic
class of the universal cover, and $\cat M= \dim M=n$ if and only if
$a^n \ne 0$.

The implication $(3) \Rightarrow (1)$ results from an obstruction
theoretic argument, \cf \cite[Lemma 8.5]{Ba}, while the implication
$(1)\Rightarrow (3)$ is obvious.
\end{proof}

\begin{theorem}[\cite{Gr1}]
\label{grom}
Every essential Riemannian manifold $M$ satisfies the inequality
\begin{equation}
\label{31}
\pisys_1(M)^n \leq C_n \vol_n(M).
\end{equation}
\qed
\end{theorem}

Here the constant $C_n$ is on the order of $n^{2n^2}$. In other words,
the quotient
\begin{equation}
\label{32}
\frac{\vol_n}{(\pisys_1)^n} > 0
\end{equation}
is bounded away from zero.  In our terminology, this implies that
homotopy systolic category is the maximal possible:
\begin{equation}
\label{33}
\syscat(M)=n \mbox{\ for\ essential\ $M$\ of\ dimension\
$n$.}
\end{equation}
Better bounds are available as the topological complexity of $M$
increases \cite{KS2}, but they are irrelevant here for the moment.
\begin{rem} 
In the appendix to \cite{Gr1}, Gromov proved an inequality of
type~\eqref{31} for Riemannian $n$-dimensional polyhedra that satisfy
condition (3) of \theoref{equiv}.
\end{rem}

The converse was proved in \cite{Ba}.  Namely, systolic category is
less than $n$ for inessential $n$-manifolds, \cf Section \ref{four}.
Combined with \theoref{equiv}, this yields the following theorem, in
harmony with equality~\eqref{amazingconjecture}.

\begin{theorem}
\label{33b}
Let $M$ be a $n$-manifold.  Then $\syscat(M)=n$ if and only if
$\cat(M)= n$.  \qed
\end{theorem}

\section{Inessential manifolds and pullback metrics}
\label{four}

We prove a compression result valid for any of our systolic notions,
including the homotopy 1-systole.

\begin{theorem}[Compression theorem]
\label{comp}
Let $n=k_1 + \ldots + k_d$.  Let $f: M^n \to K$ be map to a polyhedron
such that
\begin{enumerate}
\item
$f(M)\subset K^{(n-1)}$;
\item
$\sys_{k_i}(K)>0$ for each $k_i<n$, $i=1,\ldots, d$;
\item
the induced homomorphism $f_*$ is injective in $k_i$-dimensional
homology for each $k_i<n$, $i=1,\ldots, d$ (or in $\pi_1$ if we are
working with $\pisys_1$).
\end{enumerate}
Then the inequality 
\begin{equation}
\label{51}
\prod_i\sys_{k_i} \leq C\vol_n
\end{equation}
is violated, for any $C<\infty$, by a suitable metric on $M$.
\end{theorem}

\begin{proof}
Without loss of generality we can assume that $M$ and $K$ are
polyhedra and $f$ is linear on each simplex, and on each top
dimensional simplex $f$ is a projection to a face of positive
codimension.  Choose a fixed piecewise linear metric $\gmetric_{PL}$
on $K^{(n-1)}$.  This amounts to a choice of a positive quadratic form
on each simplex, such that the forms of neighboring simplices agree on
their common face.  Then the pullback~$f^*(\gmetric_{PL})$ is a
positive quadratic form on $M$ which is of rank at most~$(n-1)$ at
every point of $M$.  Nonetheless, volume of a Lipschitz simplex can be
defined as usual with respect to this form.  Thus its $k$-systole can
be defined.  Since $f$ induces a monomorphism $f_*$ by hypothesis, the
systole is positive by construction:
\[
\sys_k \left( M, f^*(\gmetric_{PL}) \right) \geq \sys_k \left(
K^{(n-1)}, \gmetric_{PL} \right) > 0.
\]
Then the inequality~\eqref{51} is certainly violated, to the extent
that the left hand side is positive, whereas the right hand side
vanishes.  The ``metric'' $f^*(\gmetric_{PL})$ has two shortcomings.
First, it is only defined on the full tangent space at $x\in M$ when
$x\not \in M^{(n-1)}$, \ie outside the codimension~1 skeleton of $M$.
Furthermore, the positive quadratic form is not definite, even in the
interior of a top dimensional cell.  Both shortcomings are overcome by
the metric
\begin{equation}
\label{34}
\phi_{M^{(n-1)}}^{\phantom{I}} f^*(\gmetric_{PL}) + \left(1-
\phi_{M^{(n-1)}}^{\phantom{I}} \right) \gmetric_0,
\end{equation}
where $\gmetric_0$ is a fixed smooth ``background'' metric on $M$,
while $\phi_{M^{(n-1)}}^{\phantom{I}}$ is a function of ``bump'' type,
which equals 1 outside of an $\eps$-neighborhood of the
codimension 1 skeleton, and vanishes in an $\eps/2$-neighborhood.  The
volume becomes arbitrarily small when $\eps$ tends to zero, thus
violating \eqref{51}, \cf \cite{Ba}, \cite[section~6]{BabK}.
\end{proof}

The following corollary, which is a converse of Theorem~\ref{grom}, is
due to I.~Babenko~\cite{Ba}.

\begin{corollary}
If $M$ is not essential, then it fails to satisfy Gromov's inequality
\eqref{31} for the homotopy $1$-systole.
\end{corollary}

\begin{proof}
Indeed, if $M^n$ is inessential then, by \theoref{equiv}, there exists
a map $f: M \to K=K(\pi_1 (M),1)$ that induces an isomorphism of
fundamental groups and such that $f(M) \subset K^{(n-1)}$, and we
apply our compression theorem to the homotopy 1-systole.
\end{proof}

\section{Manifolds of dimension 3}

The main result of the paper \cite{GG}, \cf also \cite{OR}, is the
following theorem.
\begin{theorem}\label{gaga}
Let $M$ be a $3$-manifold.  Then $\cat M=1$ if $M$ is simply
connected, $\cat M=2$ if $\pi_1(M)$ is a free non-trivial group, and
$\cat M=3$ otherwise.  \qed
\end{theorem}

The proof, in the orientable case, goes roughly as follows.  Decompose
a 3-manifold $M$ as a connected sum.  If $\pi_1(M)$ is not free, then
at least one of the summands has finite fundamental group or is
aspherical~\cite{H}.  Let $X$ be such a summand.  If $X$ is
aspherical, then $M$ is essential.  If~$\pi_1(X)=\pi$ is finite, we
have $\pi_2(X)=0$, since $X$ is covered by the homotopy 3-sphere.
Hence we have the Hopf exact sequence
\[
\pi_3(X) \to H_3(X) \to H_3(\pi) \to 0,
\]
and the degree of the first (Hurewicz) map is equal to the order of
fundamental group. 

In the case of free fundamental group, the theorem follows from the
fact that the manifold is homotopy equivalent to a connected sum of
2-sphere bundles over the circle.

\begin{cory}
Let $M$ be a $3$-manifold. Then: \par {\rm (i)} if $M$ is simply
connected then $\cat M=\syscat M=1$; \par {\rm (ii)} if $\pi_1(M)$ is
not a free group then $\cat M=\syscat M=3$; \par {\rm (iii)} if $M$ is
orientable and $\pi_1(M)$ is a free non-trivial group, then $\cat
M=\syscat M=2$.
\end{cory} 

For non-orientable 3-manifold with free fundamental group we only
prove that $1\le\syscat M \le \cat M=2$.  In fact, it turns out that
in this case we also have $\syscat M=2$ \cite{KR2}.
    
\begin{proof}
(i) is obvious.

\par (ii) The equality $\cat M=3$ follows from \theoref{gaga}, and we
apply \theoref{33b}.

\par (iii) If the fundamental group of $M$ is free, then $\cat M=2$ by
\theoref{gaga}.  Since $b_1(M)\geq 1$, we have $\syscat M \ge 2$ by
\eqref{gromov}, and~$\syscat M \le 2$ by \theoref{33b}.
\end{proof}

\begin{remark}
Let $f: M \to K(\pi_1(M),1)$ be a map that induces an isomorphism of
fundamental groups. An interesting role in dimension~3 is played by
the lowest dimension $d_{\min}$ of the skeleton of $K(\pi_1(M),1)$
that~$f$ can be deformed into.  For 3-manifolds, this dimension
determines both categories.  Namely, if $d_{\min}= 1$ then both
categories are~2, whereas if~$d_{\min}=3$, both categories are also 3.
Here $d_{\min}$ cannot be equal to~2. (Indeed, if $f$ can be deformed
into the 2-skeleton then~$M$ is inessential, and so $\cat M \le 2$,
and so $\pi_1(M)$ is free, and so $K(\pi_1(M),1)$ is homotopy
equivalent to a wedge of circles.) The case of maximal
$d_{\min}=n=\dim(M)$ also characterizes the case when both categories
equal $n$, by Theorem~\ref{33b}.  Could one say something in general
when this dimension is smaller than $n$?
\end{remark}

\section{Manifolds with category smaller than dimension}

Assume $b_1(M)=b \geq 1$.  An optimal systolic inequality \eqref{75}
was studied in \cite{IK}, \cf \cite{BCIK2, BCIK1, KL}.  Denote by
$J_1(M)$ the $b$-dimensional (Jacobi) torus of~$M$.  Consider the
homomorphism $\gf: \pi_1(M) \to H_1(M;\Z)\to \Z^b$, where the first
homomorphism is the abelianisation and the second one is the quotient
homomorphism over the torsion subgroup, and let
\begin{equation}
\label{81}
\AJ=\AJ_M : M \to J_1(M)
\end{equation}
be a map that induces the homomorphism $\gf$ on fundamental groups
(the so-called Abel-Jacobi map, \cf \cite{Li}).  We have the pull-back
diagram
\begin{equation}
\CD
\overline M @>\overline \AJ >> \R^b\\
@VVV @VVpV\\
M @>\AJ>> J_1(M)
\endCD
\end{equation}
where the map $p$ is the universal cover. A typical fiber $\overline
\FF$ of~$\overline{\AJ}$ projects diffeomorphically to a typical fiber
$\FF$ of $\AJ$. Denote by
\[
\left[\overline \FF \right]\in H_{n-b}(\overline M;G)
\]
the homology class of $\overline \FF$ in $\overline M$ where $G=\Z$ if
$M$ is orientable and $G=\Z_2$ if $\overline M$ is non-orientable).

Let $\gmetric$ be a metric on $M$.  Following Gromov \cite{Gr1},
denote by~$\deg(\AJ)$ the least $\gmetric$-area of an $(n-b)$-cycle
representing the homology class $[\overline \FF]$.  Note that this
notion of degree is not a topological invariant, unless of
course~$n=b$.

\begin{theorem}
\label{71z}
Let $M$ be an $n$-manifold with $b_1(M)=b$.  Then the following lower
bound for the systolic category holds:
\begin{equation}
\label{71y}
\mbox{\ if \ } [\overline \FF]\not= 0 \mbox{\ \ then\ } \syscat(M)
\geq b+1 ,
\end{equation}
where we allow a systole of $(\overline M, \overline \gmetric)$ to
participate in the definition of systolic category \eqref{syscat}.
\end{theorem}

\begin{proof}
The Jacobi torus is equipped with the stable norm~$\|\;\|$ associated
with~$\gmetric$.  Equip the Jacobi torus with the Euclidean norm
$\|\;\|_E$ defined by the ellipsoid of largest volume inscribed in the
unit ball of the stable norm $\|\;\|$.  It was proved in \cite{IK}
that the homotopy class of $\AJ$ contains a map
\begin{equation}
\label{???}
f: (M,\gmetric) \to (J_1(M), \|\;\|_E),
\end{equation}
where the map $f$ may not be distance decreasing, but is nonetheless
non-increasing on $b$-dimensional areas.  It follows from the coarea
formula that $M$ satisfies an optimal inequality
\begin{equation}
\label{75}
\stsys_1(\gmetric) ^b \deg(\AJ) \leq (\gamma_b)^{b/2}
\vol_n(\gmetric),
\end{equation}
where $\gamma_b$ is the Hermite constant, \ie the maximum, over all
lattices of unit covolume in $\R^n$, of the least square-length of a
nonzero vector in the lattice.  By hypothesis, the class of the fiber
is nontrivial, and therefore~$\sysh_{n-b}(\overline M, \overline
\gmetric) \leq \deg(\AJ)$.  It follows that its systolic category is
at least $b+1$.
\end{proof}

\begin{remark}
If $[\FF]\not=0$, then, by Poincar\'e duality, we
have 
$$
\cat M\geq \cuplength(M)\geq b+1,
$$ 
\cf \eqref{85}, in harmony with~\eqref{amazingconjecture}.  Here we
must use Poincar\'e duality with $\Z_2$-coefficients in the
non-orientable case.  Notice, however, that the condition $[\overline
\FF]\not = 0$ is more general than~$[\FF]\not = 0$, and in particular
can be satisfied even if the cup product is trivial on $M$, as in the
case of the NIL geometry on a compact 3-manifold.
\end{remark}

\begin{question}
Does a manifold with $[\overline \FF]\not = 0$ satisfy
\begin{equation}
\label{73}
\cat(M)\geq b_1(M)+1 \quad ?
\end{equation}
\end{question}

\begin{example}
\label{112}
Let $\Sigma$ be a surface different from $S^2$.  Let $M=\Sigma \times
S^n$.  Then the cup-length of $M$ is 3 for a suitable choice of
coefficients.  Hence Lusternik-Schnirelmann category is 3, and one can
ask if the systolic category of $\Sigma \times S^n$ is equal 3.  This
is true in the orientable case since real cup-length is a lower bound
for systolic category \eqref{85}.  Now consider $M= \rp^2 \times S^n$.
We conjecture that $M$ satisfies an inequality of type
\begin{equation}
\label{71}
\sys_1 \sys_1 \sys_n < C(M) \vol_{n+2},
\end{equation}
where perhaps $\sys_n$ should be replaced by either the stable
systole, or systole with $\Z_2$ coefficients.  See also \eqref{86}.
How does this fit in with the results of M.~Freedman \cite{Fr} ?
\end{example}

\section{Category of simply connected manifolds}

We first recall the following general result.
\begin{theorem}[\mbox{\cite[Theorem~1.50]{CLOT}}]
\label{condim}
For every $(k-1)$-connected $CW$-space $X$ we have $\cat X\le\dim
X/k$.  \qed
\end{theorem}

\begin{cory}\label{small}
If $M^n, n=4,5$ is a simply connected manifold which is not homotopy
equivalent to the sphere $S^n$, then $\cat M=2$.
\end{cory}

\begin{proof}
Since $M$ is not a homotopy sphere, then $H^2(M;\Z_p)\ne 0$ for some
prime $p$. Indeed, otherwise, by Poincar\'e duality, $H^{n-2}(M;\Z_p)
=0$ for all prime $p$. In other words, $H^i(M;\Z_p)=0$ for all primes
$p$ and all $i\ne 0,n$.  By the Universal Coefficient Theorem, we have
$H_i(M)=0$ for $i\ne 0,n$, and so $M$ is a homotopy sphere because it
is simply connected.

Let $x\in H^2(M;\Z_p)\setminus \{0\}$.  By Poincar\'e duality, there
exists an element~$y\in H^{n-2}(M;\Z_p)$ such that $xy\ne 0$. Now the
cup-length estimate implies $\cat M\ge 2$, and the claim follows from
\theoref{condim}.
\end{proof}

\begin{example}[Four-manifolds with category 2]
If $M^4$ is simply connected and is not homotopy equivalent to $S^4$,
then the second Betti number is positive, and $M$ satisfies the stable
systolic inequality in middle dimension:
\begin{equation}
\stsys_2(M) ^2 \leq C(b_2(M)) \vol_4(M),
\end{equation}
\cf inequalities~\eqref{gromov} and \eqref{61}, and hence
$\syscat(M)=2$.  On the other hand, $\cat M=2$ by \coryref{small}, in
harmony with equally \eqref{amazingconjecture}.
\end{example}

\begin{example}
\label{smale}
By~\cite{Sm}, a simply connected spin 5-manifold which is a rational
homology sphere is a connected sum of manifolds $M_k$, $k\geq 2$,
where
\[
H_2(M_k)=\Z_k + \Z_k.
\]
By \coryref{small}, all these manifolds have Lusternik-Schnirelmann
category equal to 2.  It is unknown whether any systolic inequalities
are satisfied by $M_k$.  Based on the value of its $\Z_k$-cup-length,
we could conjecture the existence of a systolic inequality
\begin{equation}
\label{101}
\sys_2 \sys_3 \leq C(M_k) \vol_5,
\end{equation}
where the systoles are calculated over the cycles with $\Z_k$
coefficients which are not zero-homologous, \cf \eqref{86}.  An
interesting family of metrics on $M_k$ was studied in \cite
[Theorem~8.2] {PP}.  Does it satisfy inequality~\eqref{101}?
\end{example}

\begin{example}
\label{singhof}
There are examples of $S^2$-bundles $M$ over spheres satisfying
$\cat(M)=3$ \cite[\S 3]{S}.  Such a bundle $M^{16}\to S^{14}$ can be
induced by an element of $\pi_{14}(S^4)$ from the bundle $\cp ^3 \to
S^4$ obtained as a circle quotient of the quaternionic Hopf bundle
$S^7 \to S^4$, \cf \cite{Iw}.  Let $u\in H^2(M^{16}, \Z)=\Z$ be a
generator, and let 
\[
f: M \to K(\Z, 2)=\cp ^2{\infty}
\]
classify $u$. Because of the construction of $M^{16}$, the map $f$
factors through a skeleton $\cp ^3\subset K(\Z, 2)$.

Hence $M^{16}$ admits a positive quadratic form of zero ``volume'' (as
well as zero 14-systole), but positive 2-systole, \cf
Section~\ref{four}.  By compression theorem~\ref{comp}, there is no
inequality of type $\sys_2^8 \leq C \vol_{16}$.  By Poincar\'e duality
and \eqref{52}, the manifold $M^{16}$ has stable systolic category 2.
Since there is no torsion in homology, no higher value could be
expected for systolic category, in contrast with~\eqref
{amazingconjecture}.
\end{example}

\section{Gromov's calculation for stable systoles}

In this section, we present Gromov's proof of the optimal stable
systolic inequality~\eqref{61} for complex projective space $\C P^n$,
\cf \cite[Theorem 4.36]{Gr3}, based on the cup product decomposition
of its fundamental class, \cf \eqref{86}.  In Section~\ref{massey} we
will adapt this calculation to a situation where cup product is
trivial.  Following Gromov's notation, let~$\alpha\in H_2(\cp^n;\Z)$
be a generator in homology, and $\omega\in H^2(\cp^n;\Z)$ the dual
generator in cohomology.  Then $\omega^n$ is a generator of
$H^{2n}(\cp^n,\Z)$.  Let~$\eta \in \omega$ be a closed differential
2-form.  Then
\begin{equation}
\label{61b}
1= \int_{\cp^n} \eta^n .
\end{equation}
Now let $\gmetric$ be a metric on $\cp^n$.  Then
\begin{equation}
\label{61d}
1 \leq n!  \left( \| \eta \|_\infty \right)^n \vol_{2n}(\gmetric),
\end{equation}
where $\| \; \|_\infty$ is the comass norm on forms (see \cite{Gr3}
for a discussion of the constant).  Here the comass norm of a
differential $k$-form is the supremum of pointwise comass norms.  The
pointwise comass norm for decomposable linear $k$-forms coincides with
the natural Euclidean norm on $k$-forms associated with $\gmetric$.  In
general it can be defined by evaluating on $k$-tuples of unit vectors
and taking the maximal value, \cf \cite{Fe, BanK}.  Therefore
\begin{equation}
\label{61c}
1 \leq n! \left( \| \omega \|^* \right)^n \vol_{2n}(\gmetric),
\end{equation}
where $\|\;\|^*$ is the comass norm in cohomology, obtained by
minimizing the norm~$\|\;\|_{\infty}$ over all closed forms
representing the given cohomology class.  Denote by~$\|\;\|$ the
stable norm in homology.  Recall that the normed lattices $(H_2(M,
\Z), \|\;\|)$ and $(H^2(M, \Z), \|\;\|^*)$ are dual to each other
\cite{Fe, BanK}.  Therefore
\[
\|\alpha \| = \frac{1}{\|\omega \|^*},
\]
and hence
\begin{equation}
\label{61}
\stsys_2(\gmetric)^n = \| \alpha \|^n \leq n! \vol_{2n}(\gmetric) .
\end{equation}
Moreover the inequality so obtained is sharp (equality is attained by
the 2-point homogeneous Fubini-Study metric).  In our terminology,
this calculation can be summarized by writing
\[
\syscat(\cp^n) = n,
\]
which agrees with $\cat(\cp^n)$, in harmony with
equality~\eqref{amazingconjecture}.

\m
Gromov also proved a stable systolic inequality associated with any
decomposition of the real fundamental cohomology class of $M$ as a cup
product, \cf \cite{Gr1, BanK}.  This can be restated in terms of
category as follows:
\begin{equation}
\label{85}
\syscat(M) \geq \cuplength_\R(M).
\end{equation}
Pu's inequality \eqref{21} for $\rp^2$ as well as its generalisations
by Gromov lend support to the following conjecture:
\begin{equation}
\label{86}
\syscat(M) \geq \cuplength (M),
\end{equation}  
(for the appropriate choice of systoles), which would be implied by a
corresponding equality in~\eqref{amazingconjecture}.  Interesting test
cases are \eqref{71}, \eqref{101}.

\section{Massey products via DGA algebras}

Let $(A,d)$ be a differential graded associative (=DGA) algebra with a
differential $d$ of degree 1, and let $H=\Ker d/\IM d$ be the homology
algebra of $(A,d)$. We assume that the induced product in the graded
algebra $H$ is skew-commutative. Given (homogeneous) $u,v, w\in H$
with $uv=0=vw$, the triple Massey product
\[
\la u,v,w\ra = \la u,v,w\ra_A \subset H
\]
is defined as follows.  Let $a,b,c$ be elements of $A$ whose
homology classes are $u,v,w$ respectively. Then $dx=ab$, $dy=bc$ for
suitable $x,y\in A$, and $\la u,v,w\ra$ is the set of elements of the
form
\[
xc - (-1)^{|u|}ay,
\]
see \cite{M, RT} for more details.  The set $\la u,v,w\ra$ is a coset
with respect to the so-called {\it indeterminacy subgroup} $\In
\subset H^{|u|+|v|+|w|-1}$,
\[
\In = uH^{|v|+|w|-1}+H^{|u|+|v|-1|}w.
\]

\begin{prop}\label{real}
Let $X$ be a finite $CW$-space.  Let $C^*(X;\Z)$ and $C^*(X;\R)$ be
the singular cochain complexes with coefficients, respectively, in
$\Z$ and~$\R$.  We equip these chain complexes with the
Alexander--Whitney product.  Given $u,v,w\in H^*(X;\Z)$, suppose that
the Massey product $\la u,v,w\ra\subset H^*(M;\Z)$ is defined and let
$\In$ be its indeterminacy subgroup.  Then the Massey product
\[
\la u_{\R},v_{\R},w_{\R}\ra \subset H^*(M;\R)
\]
is defined, and we have 
\[
\la u_{\R},v_{\R},w_{\R}\ra =\la u,v,w\ra_{\R}+\In \otimes \R.
\]
\end{prop}

\begin{proof} Clearly, if $x\in\la u,v,w\ra$ then $x_{\R}\in\la 
u_{\R},v_{\R},w_{\R}\ra$. Furthermore, the indeterminacy subgroup for $\la 
u_{\R},v_{\R},w_{\R}\ra$ is the subgroup $$
u_{\R}H^{|v|+|w|-1}(X;\R)+H^{|u|+|v|-1|}(X;\R)w_{\R},
$$ 
and the result follows.
\end{proof}

\m It is clear that if $f: A\to A'$ is a morphism of DGA algebras such
that $f_*: H \to H'$ is an isomorphism, then $f_*: H \to H'$ induces
an isomorphism of Massey products. However, one can relax the
requirement for $f$ to be a ring homomorphism.

\begin{df}
[\cf \cite{M, BG}] We say that an additive chain map $f: A \to A'$ is
{\it a ring map up to higher homotopies\/} if there exists a family of
linear maps $f_i: A^{\otimes i} \to A'$ of degree $1-i$, where $f_0=f$
and
\begin{equation}\label{higher}
\begin{aligned}
d\circ f_i & +(-1)^if_i\circ d= \cr & = \sum _{j=1}^{i-1} (-1)^j
(\mu(f_j \otimes f_{i-j}-f_{i-1}(1^{j-1}\otimes \mu\otimes 1^{i-j-1}))
\end{aligned}
\end{equation}
Here $d$ and $\mu$ are the differential and the multiplication in $A'$, 
respectively, and $1^{\otimes k}$ is the identity map of $A^{\otimes k}$.
\end{df}

\begin{theorem}\label{may}
 If $f: A \to A'$ is a ring map up to higher homotopies, then the map
$f_*: H \to H'$ is a ring homomorphism. Moreover, if a Massey triple
product $\la u,v,w\ra$ is defined for some $u,v,w\in H$, then the
Massey triple product $\la f_*u,f_*v, f_*w\ra$ is defined, and
\begin{equation}\label{induced} 
f_*(\la u,v,w\ra)\subset\la f_*u,f_*v, f_*w\ra. 
\end{equation}
Finally, if $f_*: H \to H'$ is an isomorphism then  
$$
f_*(\la u,v,w\ra) = \la f_*u,f_*v, f_*w\ra.
$$
\end{theorem}

\begin{proof}
Since $df_1-f_1d=0$ and $df_2+f_2d=f_1\mu -\mu(f_1\otimes f_1)$, we
conclude that $f$ induces a ring homomorphism $H \to H'$. The claims
on Massey products can be proved directly but a bit tediously, \cf
\cite[Theorem~1.5]{M} where more general results are proved.
\end{proof}

\begin{rem} 
As it is clear from the proof, the existence of $f_0, f_1, f_2$
already implies that $f_*$ is a ring homomorphism. Moreover, the
existence of $f_0, f_1, f_2$ yields formula \eqref{induced}.  In fact,
one can $n$-tuple Massey products, and an analog of \eqref{induced}
holds provided the maps $f_1, \ldots, f_n$ as in \eqref{higher} exist
\cite[Theorem 1.5]{M}.
\end{rem}

\m Given a manifold $M$, let $\Omega^*(M)$ be the de Rham algebra of
differential forms on $M$, and let $C^*(M;\R)$ be the singular cochain
with the Alexander--Whitney product.  Let
\[
\rho: \Omega^*(M) \to C^*(M;\R)=\Hom (C_*(M;\Z), \R )
\]
be the map given by integration of a form over a chain.  (To be
rigorous, we must consider chains in $C_*(M)$ generated by {\it
smooth} simplices, but one can prove that the corresponding chain
complex yields the standard homology group.) By the de~Rham Theorem,
the map $\rho$ induces an isomorphism of homology of these two
complexes, \ie an isomorphism $H^*_{\dR}(M) \to H^*(M;\R)$.

Clearly, the map $\rho$ is not a ring homomorphism (since the
algebra~$C^*(M;\R)$ is not commutative), but it is a ring map up to
higher homotopies \cite[Proposition~3.3]{BG}.  Therefore, by
\theoref{may}, the map $\rho_*:H^*_{\dR}(M) \to H^*(M;\R)$ induces a
bijection of Massey products.

Consider the map
\begin{equation}\label{r}
\CD
r:H^*(M;\Z) @>>> H^*(M;\R) @>\rho_* >> H^*_{dR}(M)
\endCD
\end{equation}
where the first map has the form
$$
H^*(M;\Z) \to H^*(M;\Z) \otimes \R=H^*(M;\R),\quad a \mapsto a\otimes 1.
$$

\begin{df}
We say that a de Rham cohomology class is {\it integral} if it belongs to the 
image of the map $r$ as in \eqref{r}.
\end{df}

\begin{lemma}\label{integral}
Let $u,v,w\in H^*_{\dR}(M)$ be integral classes and assume that the groups 
$H^{|u|+|v|}(M;\Z)$ and $H^{|v|+|w|}(M;\Z)$ are torsion free. If the Massey 
product $\la u,v,w\ra\subset H^*_{\dR}(M)$ is defined, then it contains an 
integral class.
\end{lemma}

\begin{proof} Let $\ov u, \ov v, \ov w\in H^*(M;\Z)$ be elements such that 
$r(\ov x)=x$ for~$x=u,v,w$.  By our torsion hypotheses, we conclude
that $\ov u\, \ov v=0=\ov v\,\ov w$, and therefore the Massey product $\la
\ov u,\ov v, \ov w\ra\subset H^*(M;\Z)$ is defined. But, by
\propref{real} and \theoref{may}, we have the inclusion 
\[
r( \la \ov u,\ov v, \ov w\ra)\subset \la u,v,w\ra,
\]
and the lemma follows from the fact that the set $r( \la \ov u,\ov v,
\ov w\ra)$ consists of integral classes.
\end{proof}

\section{Gromov's calculation in the presence of a Massey product} 
\label{massey}

We will follow Gromov's suggestion \cite[7.4.$C'$, p.~96]{Gr1} of
exploiting nontrivial Massey products in combination with
isoperimetric quotients, so as to obtain geometric inequalities.  Some
examples appear in \cite{DR}.  We will use in an essential way the
identification of two distinct Massey product theories, \cf
\eqref{102}.

Given a metric $\gmetric$ on $M$, let $ \IQ(M,\gmetric) $ be its
isoperimetric quotient, defined as the maximin of quotients of the
comass norm of a $(k-1)$-dimensional differential form which is a
primitive, by the comass norm of an exact $k$-form on $M$, over all
$k$.  Namely,
\begin{equation}
\IQ(\gmetric) = \max_k \; \sup_{\alpha\in \Omega^k} \inf_\beta \left\{
\left.  \frac{\|\beta\|^*} {\|\alpha\|^*} \; \right| d \beta = \alpha
\right\}.
\end{equation}

\begin{theorem}
\label{91}
Let $M$ be an orientable $n$-manifold and $p_1,p_2$ positive integers.
Let $\omega_i \in H^{p_i}_{\dR}(M), i=1,2$ be the image of an integral
generator under the map $r$, as in $\eqref{r}$.  Suppose the Massey
triple product $\langle \omega_1, \omega_1, \omega_2 \rangle$ is
defined and nontrivial.  Let 
\[
p_3=n-(2p_1+p_2-1),
\]
and assume $b_{p_1}= b_{p_2}=b_{p_3}=1$.  Assume furthermore that the
group $H^{p_1 + p_j}(M;\Z),$ $j=1,2$ is torsion free.  Set
\[
A_j=\frac {n!}  { {p_j!\; (p_1+p_j)!\;
p_3!}}  \binom {p_1 + p_j}{p_1}, \quad j=1,2 .
\]
Then every metric $\gmetric$ on $M$ satisfies the inequality
\[
 {\stsys_{p_1}(\gmetric)^2 \stsys_{p_2}(\gmetric)
\stsys_{p_3}(\gmetric)} \leq (A_1+A_2) {\IQ(\gmetric) }
\vol_n(\gmetric).
\]
\end{theorem}

\begin{proof}
Throughout the proof, we will denote the de Rham cohomology
group~$H^*_{\dR}(M)$ by $H^*(M)$.  We have $b_{n-p_3}=b_{p_3}=1$ by
Poincar\'e duality.  Recall that the indeterminacy of the Massey
product $\langle \omega_1,\omega_1, \omega_2\rangle$ is an $\R$-vector
subspace of $H^{n-p_3}(M;\R) =\R$.  Since the Massey product is
assumed nontrivial, it must have zero indeterminacy.

\m
Choose a representative $\eta_i\in \omega_i$.  By hypothesis, the
$(p_1+p_j)$-form $\eta_{1} \wedge \eta_{j}, j=1,2$ is exact.  Choose a
primitive $\eta_{1j}$ satisfying $d \eta_{1j} = \eta_{1} \wedge
\eta_j$ of minimal comass norm, so that

\begin{equation}
\label{64}
\frac {\| \eta_{1j} \|_\infty} {\| \eta_{1} \wedge \eta_j \|_\infty}
\leq \IQ(M,\gmetric) .
\end{equation}

Notice that the cohomology class $\omega$ of the closed form
\[
\eta_{11} \wedge \eta_2 -(-1)^{p_1} \eta_1 \wedge \eta_{12}
\]
is equal to the Massey product $\langle \omega_1, \omega_1, \omega_2
\rangle _{\dR} \in H^{n-p_3}(M;\R)$.  By Poincar\'e duality, there is
a class $\omega_3 \in H^{p_3}(M;\R)=\R$ which is the image of an
integral generator, such that $ (\omega_3\cup\omega) [M] >0.  $ By
\lemref{integral}, the class $\omega$ is integral, and hence the class
$\omega_3\cup\omega$ is integral, so that
\begin{equation}
\label{102}
(\omega_3\cup\omega ) [M] \ge 1.
\end{equation}
Choose a $p_3$-form $\eta_3$ in the cohomology class $\omega_3$.  Then
\begin{equation}
1 \leq \int _{M} \left( \eta_{11} \wedge \eta_2 -
{(-1)}_{\phantom{I_{}}}^{p_1} \eta_1 \wedge \eta_{12} \right) \wedge
\eta_3 .
\end{equation}
Similarly to \eqref{61d} and exploiting \eqref{64}, we obtain
\[
\begin{aligned}
1 & \leq \sum_{j=1}^2 \frac { n!}{ {p_j!\; (p_1+p_j-1)!\; p_3!}}  \|
\eta_j \|_\infty \; \| {\eta_{1j}} \|_\infty \; \| \eta_3 \|_\infty
\vol_{n}(\gmetric) \cr & \leq (A_1 + A_2) \IQ(\gmetric) \; \| \eta_1
\|_\infty \; \| {\eta_1} \|_\infty \; \| \eta_2 \|_\infty \| \; \eta_3
\|_\infty \vol_n(\gmetric).
\end{aligned}
\]
Hence, similarly to \eqref{61c}, we obtain
\[
1 \leq (A_1+A_2) \IQ(\gmetric) \left( \| \omega_1 \|^*\right)^2 \; \|
{\omega_2} \|^* \; \| \omega_3 \|^* \vol_n(\gmetric) .
\]
Now let $\alpha_i \in H_{p_i}(M;\R)$ be the generator dual to
$\omega_i$.  By duality of comass and stable norm, and since
$b_{p_i}=1$, we have $\| \alpha_i\| =\frac{1}{ \| \omega_i
\|^*}$. Thus, similarly to \eqref{61} we obtain the inequality
\[
\stsys_{p_1}(\gmetric)^2 \stsys_{p_2}(\gmetric) \stsys_{p_3}(\gmetric)
\leq (A_1+A_2)\IQ(\gmetric) \; \vol_n(\gmetric) ,
\]
as required.
\end{proof}

\section{A homogeneous example}
Consider the homogeneous manifold
\[
M^{n^2+2n-16}= SU(n+1)/SU(3)\times SU(3),
\]
where $n\geq 5$ \cite{Tr}, \cite[Chapter~11]{GHV}.  Denote
by~$\omega_i$ the generator of $H^i(M)$, for $i=4,6$.  The Massey
triple product
\[
\langle \omega_4, \omega_4, \omega_6 \rangle \in H^{13}(M,\Z)
\]
is nontrivial.  For $n=5$, we obtain that every metric $\gmetric$ on
the manifold $SU(6)/SU(3)\times SU(3)$ satisfies
\begin{equation}
\label{92}
\frac {\stsys_4(\gmetric)^2 \stsys_6(\gmetric)^2} {\IQ(\gmetric)} \leq
19! \vol_{19}(\gmetric) .
\end{equation}

\begin{remark}
This inequality is not exactly of the form envisioned in~\eqref{dd},
since the constant is metric-dependent via the $\IQ$ factor.  Applying
this technique to the standard 3-dimensional nilmanifold, we get a
lower bound given by a product of {\em four\/} systoles, for the
quantity $\IQ(\gmetric) \vol(\gmetric)$.  Since we expect the category
of a 3-manifold to be at most 3, it is natural to define an analog of
systolic category in the presence of a Massey product, by {\em
subtracting\/} the number of $\IQ$ factors, as in \eqref{94} below,
appearing as the result of integrations involved in defining a Massey
product.  This would be consistent with the lower bound for $\cat(M)$
in terms of weights of Massey products.  Thus, the inequality can be
restated as a lower bound for an IQ-modified systolic category:
\begin{equation}
\label{94}
{\rm cat}_{\rm sys} ^{\iqvol} (M^{19}) \geq 4-1 =3.
\end{equation}
\end{remark}
On the LS side, it can be shown that $3\le \cat M^{19}\le 4$, and we
expect that that $\cat =4$.  Meanwhile, John Oprea (private
communication) proved the following proposition.  The definition and
properties of rational Lusternik--Schnirelmann category can be found
in \cite{CLOT}.

\begin{prop}
\label{104}
The rational Lusternik--Schnirelmann category of the manifold~$M^{19}$
is equal to $3$.
\end{prop}

\begin{proof}
For every simply connected manifold~$X$, the rational
Lusternik--Schnirelmann category of $X$ is equal to its rational
Toomer invariant~$e_0(X)$, see \cite{FHL}. To compute $e_0(X)$, notice
that the Sullivan model of the space $SU(6)/(SU(3) \times SU(3))$ has
the form
$$
\Lambda(x_4,x_6,y_7,y_9,y_{11})
$$
with differential
$$
dx_i = 0,\ dy_7 = x_4^2,\ dy_9 = x_4 x_6,\ 
dy_{11} = x_6^2.
$$

Therefore a top class is given by
\begin{equation}
x_4^2 y_{11} - x_4 x_6 y_9,
\end{equation}
and so $e_0(M^{19})=3$ by \cite[Prop. 5.23]{CLOT}. 
\end{proof}

\section{Acknowledgments}
We are grateful to J. Oprea for providing a proof of Proposition~\ref
{104}.

\vfill\eject

\end{document}